\documentclass[12pt]{article}
\usepackage{amsmath,amssymb}
\usepackage[dvips]{graphicx}
\usepackage[active]{srcltx}
\usepackage[final]{epsfig}

\newcommand{\bbRP}{\mathbb{RP}}
\newcommand{\Diff}{\mathrm{Diff}}
\newtheorem{theorem}{Theorem}[section]

\newtheorem{theorem-definition}[theorem]{Theorem-Definition}
\newtheorem{theorem-construction}[theorem]{Theorem-Construction}
\newtheorem{lemma-definition}[theorem]{Lemma--Definition}
\newtheorem{proposition-definition}[theorem]{Proposition--Definition}
\newtheorem{lemma}[theorem]{Lemma}
\newtheorem{proposition}[theorem]{Proposition}
\newtheorem{corollary}[theorem]{Corollary}
\newtheorem{conjecture}[theorem]{Conjecture}
\newtheorem{definition}[theorem]{Definition}

\begin{document}
\newcommand{\Z}{{\mathbb Z}}
\newcommand{\R}{{\mathbb R}}
\newcommand{\Q}{{\mathbb Q}}
\newcommand{\C}{{\mathbb C}}
\newcommand{\lms}{\longmapsto}
\newcommand{\lra}{\longrightarrow}
\newcommand{\hra}{\hookrightarrow}
\newcommand{\ra}{\rightarrow}
\newcommand{\sgn}{\rm sgn}

\begin{titlepage}
\title{Moduli spaces of 
convex projective structures on surfaces}  
\author{V.V. Fock, A. B. Goncharov}
\end{titlepage}
\date{{\it \small To the memory of Yan Kogan}}
\maketitle
\tableofcontents
\vskip 4mm \noindent
 
\section{Introduction}

Let $S$ be an orientable compact smooth surface possibly with boundary. A projective structure on $S$ is defined by an atlas on the interior $S_0$ of $S$ whose transition functions are given by (restrictions) of projective transformations.  For any projective structure on $S$ one can associate a developing map 
$$
\varphi: \widetilde{S_0}\rightarrow \R {\Bbb P}^2
$$ 
where $\widetilde{S_0}$ is a universal cover of $S_0$. It is defined up to the right action of $PSL_3({\mathbb R})$. 

An open domain $D$ in $\R {\Bbb P}^2$ is {\it convex} if any line intersects $D$ by a connected interval, possibly empty. A projective structure on $S$ is {\em convex} if the developing map $\varphi(\widetilde{S}_0)$ is an embedding and its image is a convex domain. Since $\varphi$ is injective, $\varphi(\widetilde S_0)$ is orientable. Conversely, let $\Delta$ be a discrete subgroup of $PSL_3(\mathbb R)$ such that its action on $\mathbb{RP}^2$ restricts to a free action on a convex domain $D$.  Then the quotient $D/\Delta$ is a surface with a convex projective structure. The natural isomorphism $\mu$ from $\pi_1(S)$ to $\Delta$ is called the monodromy map of a projective structure and is defined up to conjugation.

The space of convex projective structures on smooth compact surfaces without boundary was studied by W.Goldman, S.Choi  \cite{Goldman},\cite{Choi-Goldman}, 
F.Labourie \cite{Labourie}, H. Kim \cite{Kim}, J.Loftin \cite{Loftin} and implicitly N.Hitchin \cite{Hitchin}. Goldman's paper \cite{Goldman} discusses convex projective structures with geodesic boundary and regular semi-simple holonomy for the boundary. Although the main applications in that paper concerned closed surfaces, the proof involved computing the deformation space for three holed sphere. The main result stated there holds for all compact oriented surfaces with negative Euler characteristic, with these boundary conditions. Choi later extended it to non-orientable surfaces and discussed other cases of boundary holonomies. 

On the other hand, given a split semi-simple algebraic group $G$ with trivial center and a surface $S$ with non-empty boundary, we defined  in \cite{Higher} the Higher Teichm{\"u}ller space ${\cal X}_{G, S}(\R_{>0})$.  It comes equipped with a distinguished collection of coordinate systems, parametrized by the set of the isotopy classes of trivalent graphs embedded to $S$ and homotopically equivalent to $S$, plus some extra data. The mapping class group acts in these coordinates in an explicit way. We have shown in Section 10 of {\em loc.\! cit.} that for $G= PGL_2$ we recover the classical Teichm{\"u}ller spaces. 

In this paper we investigate convex projective structures on hyperbolic surfaces with a non-empty piecewise geodesic boundary. We introduce a distinguished collection of global coordinate systems on the corresponding moduli space ${\cal T}^+_3(S)$. Each of them identifies it with $\R^{8\chi(S)}$. We show that the moduli space ${\cal T}^+_3(S)$ is naturally isomorphic to the space  ${\cal X}_{G, S}(\R_{>0})$ for $G = PGL_3$. We prove that the corresponding monodromy representations are discrete and regular hyperbolic. We introduce and study the universal higher Teichm\"uller space $\mathcal T_3$ which contains the moduli spaces $\mathcal T^+_3(S)$ for all surfaces $S$. The Thompson group $\mathbb T$ acts by its automorphisms. We show that the universal higher Teichm\"uller space $\mathcal T_3$ is the  set of $\mathbb R_{>0}$-points of  certain infinite dimensional 
cluster $\mathcal X$-variety, as defined in \cite{Cluster}. Moreover the Thompson group is a subgroup of the mapping class group $\mathbb T_3$ of this cluster variety. For the spaces $\mathcal T_3^+(S)$ the situation is similar but a bit more complicated: they have orbi-cluster structure, and the mapping class group of $S$ is a subgroup of the corresponding group for the cluster ${\cal X}$-space. We quantize (= define a non-commutative deformation of the algebra
of functions on) the above 
moduli spaces.  The Thompson group acts by automorphisms of the quantum universal 
Teichm\"uller space ${\cal T}^q_3$. The space ${\cal T}^q_3$ 
can be viewed as a combinatorial version of the $W$-algebra for $SL_3$. 

We tried to make the exposition self-contained and as elementary as possible. Therefore this paper can serve as an elementary introduction to \cite{Higher}, where the results of this paper were generalized to the case when $PGL_3$ is replaced by an arbitrary split reductive group $G$ with trivial center.    

\section{Convex projective structures with geodesic boundary}

\paragraph{1. The moduli space of framed convex projective structures with geodesic boundary on $S$ and its versions.} 
A {\it convex curve} in $\R {\Bbb P}^2$ is a curve 
such that any line intersects 
it either by a connected line segment, or in no more then two points. 
If every line intersects a curve in no more then two points, it is a 
{\it strictly convex curve} in $\R {\Bbb P}^2$.  

There is a natural bijective correspondence between convex 
curves and convex domains in $\R {\Bbb P}^2$. Indeed, let $D$ be a 
convex domain in $\R {\Bbb P}^2$. Then its boundary $\partial D$ is a 
convex curve. Conversely, let $K$ be a convex curve. Then barring the two 
trivial cases when $K$ is a line or an empty set, 
the complement to $K$ is a union of two convex domains. 
One of them is not orientable, and contains  lines. The other 
is orientable, and does not contain any  line. 
The latter is called {\it interior} of the convex curve $K$, and denoted $D_K$.  
If $K$ is a line or the empty set,  the convex domain $D_K$, by definition, is 
$\R^2$, or $\R {\Bbb P}^2$. 

Given a convex domain $D \subset \R {\Bbb P}^2$ we define the projectively
dual domain  $\widehat D \subset \R \widehat {\Bbb P}^2$ as the set of
all lines in $\R{\Bbb P}^2$ which do not intersect $D$. It is a convex domain. 
The dual to
$\widehat D$ is $D$.

A curve on a surface $S$ with projective structure is called {\it geodesic} if it is a straight line segment in any projective coordinate system. Therefore a geodesic develops into an infinite collection of line segments. They are permuted by the monodromy group $\mu(\pi_1(S))$. 

If the surface $S$ has boundary the space of projective structures on $S$ is infinite dimensional even with the convexity requirement. For example the set of projective structure on a disk coincides with the set of convex domains in $\mathbb{RP}^2$ up to the action of $PSL(2,\mathbb R)$, which is obviously infinite dimensional. Therefore we need to impose more strict boundary conditions to ensure more finite-dimensional moduli space. We require the developing map to be extendable to the boundary and the image of every boundary component to be either a segment of a line or a point.  The second case is called degenerate and the boundary is called {\em cuspidal}. We say that a boundary component of a surface with a projective structure is {\em geodesic} if it the projective structure falls into one of these two cases.

A {\em framing} a projective structure with geodesic boundary is an orientation of all non-degenerate boundary components. 

\begin{definition} $\mathcal T^+_3(S)$ is the moduli space of framed convex 
real projective structures on an oriented surface $S$ with geodesic boundary considered 
up to the action of the group $\Diff^0(S)$ of diffeomorphisms isotopic to the 
identity. 
\end{definition}
(The index $3$ in $\mathcal T^+_3(S)$ indicates  the group $PSL(3,{\mathbb R})$and $+$ stands for framing). 

The space $\mathcal T^+_3(S)$ is a $2^s:1$ cover of the space $\mathcal T_3(S)$ of non-framed convex real projective structures on $S$ with geodesic boundary, where $s$ is the number of the holes,  ramified over the surfaces which have at least one cuspidal boundary component. On the other hand the orientation of $S$ provides orientations of the boundary components. Thus there is a canonical embedding  $\mathcal T^+_3(S) \hra \mathcal T_3(S)$ as a subspace with corners. 

The spaces  $\mathcal T^+_3(S)$ enjoys the following structures and properties:

{\bf 1.} There is an embedding $i: \mathcal T_2^+(S) \hra \mathcal T^+_3(S)$, where  $\mathcal T_2^+(S)$ is the classical Teichm{\"u}ller space parametrising complex structures on $S_0$ modulo the action of $\Diff^0(S)$, with chosen orientation of non-degenerate boundary components. Indeed, due to the Poincar{\'e} uniformisation theorem any Riemann surface can be represented as a quotient of the hyperbolic plane by a discrete group. Consider the Klein model of the hyperbolic plane in the interior of a conic in $\mathbb RP^2$. The geodesics are straight lines in this model. The quotient inherits the real projective structure. The orientations of the boundary components are inherited trivially.

{\bf 2.}
The map $\mu$ from $\mathcal T^+_3(S)$ to the space ${\mathcal R}_3(S)$ of homomorphisms $\pi_1(S) \to PSL_3({\mathbb R})$ considered up to a conjugation. The space $\mathcal R_3(S)$ possesses a Poisson structure (\cite{Goldman}, \cite{Goldman2}, \cite{FR}). Since $\mu$ is a local diffeomorphism, it induces a Poisson structure on $\mathcal T^+_3(S)$. 

{\bf 3.}
The involution $\sigma: \mathcal T^+_3(S)\rightarrow \mathcal T^+_3(S)$, defined by the property that the convex domain corresponding to the point $\sigma x$, where $x \in \mathcal T^+_3(S)$, is projectively dual to the convex domain corresponding to  $x$. The representation $\mu(\sigma(x))$ is defined by composing $\mu(x)$ with  the outer automorphism $g \to (g^t)^{-1}$ of the group $PSL_3(\mathbb R)$. The map $\sigma$ preserves the Poisson structure. 

{\bf 4.} 
The action of the mapping class group $\Gamma_S:= \Diff(S)/\Diff^0(S)$ on $\mathcal T^+_3(S)$. It also preserves the Poisson structure.


Before we proceed any further, let us introduce a toy model of the moduli space $\mathcal T^+_3(S)$, which does not only contains the main features of the latter, but also plays a key role in its study.  
\vskip 3mm
\paragraph{2. The moduli spaces of pairs of convex 
polygons, one inscribed into the other.} 
Let $\mathcal P_3^n$ be the space of pairs of convex  $n$-gons in $\mathbb{RP}^2$, one  inscribed into the other, and considered up to the action of $PSL_3(\mathbb R)$. 

One can think of this space  as of a kind of discrete approximation to the space of parametrised closed convex curves in $\mathbb{RP}^2$. Indeed, fix a set $R$ of $n$ points on the standard circle $S^1$.  Then for any convex curve $\gamma:S^1 \rightarrow \mathbb{RP}^2$ one can associate the convex polygon with vertices $\gamma(R)$ and the polygon with edges tangent to $\gamma$ at $\gamma(R)$.

This space has a natural Poisson structure,  and there are analogs of the maps
$\mu$,$\sigma$, $i$ and the mapping class group action. 

Namely, let $F_3$ be the space of 
flags in ${\mathbb R}^3$, or, equivalently the space of pairs 
$(A,a)$, where $a\subset\bbRP^2$ is a line, and $A\in a$ is a point. 
There exists a natural map $\mu:{\cal P}_3^n \rightarrow F^n_3/PSL_3({\mathbb R})$ 
which is the analog of the map $\mu$ for ${\cal T}_3(S)$. Its image is a 
connected component in the space of collections of flags in general position. 

The explicit formulae given below provide a definition 
of the Poisson structure in this case. 

The projective duality interchanges the inscribed and circumscribed polygons. 
It acts as an involution $\sigma$ of ${\cal P}_3^n$. 

Let ${\cal P}_2^n$ be the 
configuration space of $n$ ordered 
points on the oriented real projective line, such that the  induced cyclic 
order is compatible with the chosen orientation of $\bbRP^1$. 
Then there is a canonical embedding $i:{\cal P}_2^n \hra {\cal P}_3^n$, which 
identifies ${\cal P}_2^n$ with the space of polygons inscribed into a conic.
 Indeed, for each such a polygon we assign the circumscribed polygon given 
by tangents to the conic at the vertices of the original polygon. 
The set of stable points of $\sigma$ is precisely the image of $i$.

The role of the mapping class group is
 played by the cyclic group shifting simultaneously the vertices 
of the two polygons. Namely, we assume that the vertices of the two polygons are 
ordered, so that the induced cyclic orders are compatible with the natural orientations of the polygons. One may assume in addition that the first vertex of the 
inscribed polygon is inside of the first side of the circumscribed one. The 
generator of the cyclic group shifts cyclically the order by one.

One can unify these two moduli spaces by considering the moduli space $\mathcal T_3(\widehat S)$ of framed convex projective structures on surfaces with geodesic boundary, equipped with a finite (possible empty) collection of marked points on the boundary, see Section 2.7.  

Our goal is to introduce  a set of global coordinates on ${\cal T}^+_3(S)$, describe its natural Poisson structure in terms of these parametrization and give explicit formulae for the maps $i,\mu,\sigma$ and the action of the mapping class group. Before addressing this problem, let us first solve this problem for the toy model $\mathcal P_3^n$ since it contains most of the tools used for the problem in question.
\vskip 3mm
\paragraph{3. Parameterizations of the spaces ${\cal P}_3^n$.} 
Cut the inscribed polygons into triangles and mark two distinct points on every edge of the triangulation except the edges of the polygon. Mark also one point inside each triangle. 

\begin{theorem} \label{1.20.04.5}
There exists a canonical bijective correspondence between the space ${\cal P}_3^n$ and assignments of positive real numbers to the marked points.
\end{theorem}

{\bf Proof}. It is  constructive: we are going to describe how to construct numbers 
from a pair of polygons and visa versa.  

\begin{figure}[ht]
\centerline{\epsfbox{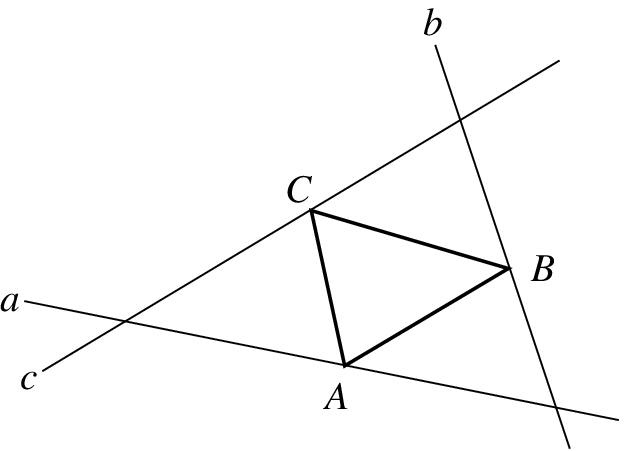}}
\caption{${\cal P}_3^3$}
\label{triangle}
\end{figure}

A small remark about notations. We denote points (resp. lines) by uppercase (resp. lowercase) letters.  A triangle on $\bbRP^2$ is determined neither by vertices nor by sides, since there exists four triangles for any generic triple of vertices or sides. If the triangle is shown on a figure, it is clear which one corresponds to the vertices since only one of the four fits entirely into the drawing. If we want to indicate a triangle which does not fit, we add in braces a point which belongs to the interior of the triangle. For example, the points $A,B$ and $C$ on Figure  \ref{triangle} are vertices of the triangles $ABC$, $ABC\{a\cap c\}$, $ABC\{a\cap b\}$ and $ABC\{b\cap c\}$.

Another convention concerns the cross-ratio. We assume  that the cross-ratio of 
four points $x_1, x_2, x_3, x_4$ on  a line is the value at $x_4$ of a projective coordinate 
taking value $\infty$ at $x_1$, $-1$ at $x_2$, and $0$ at $x_3$. So we employ the formula 
$\frac{(x_1-x_2)(x_3-x_4)}{(x_1-x_4)(x_2-x_3)}$ 
for the cross-ratio. Observe that the 
lines passing through a point form a projective line. So the cross-ratio of 
an ordered quadruple of lines passing through a point is defined. 

Let us first consider the case of ${\cal P}_3^3$. It is the space of pairs of triangles $abc$ and $ABC$ (Figure  \ref{triangle}), where the second is inscribed into the first and considered up to projective transformations. This space is one dimensional and its invariant $X$ is just the cross-ratio of the quadruple of lines $a,AB,A(b\cap c),AC$. Such invariant is called the {\em triple ratio} (\cite{Gon1}, Section 3) and is defined for any generic triple of flags in $\mathbb{RP}^2$ and can be also defined as follows. Let $\mathbb R^3$ be the three dimensional real vector space whose projectivisation is $\R {\Bbb P}^2$. Choose linear functionals $f_a, f_b, f_c \in (\mathbb R^3)^*$ defining the lines $a, b, c$. Choose non-zero vectors $\widetilde A, \widetilde B, \widetilde C$ projecting onto the points 
$A, B, C$. Then 
$$
X :=  \frac{f_a(\widetilde B)f_b(\widetilde C)f_c(\widetilde A)}{f_a(\widetilde C)f_b(\widetilde A)f_c(\widetilde B)}
$$
It is obviously independent on the choices involved in the definition, and manifestly ${\mathbb Z}/3{\mathbb Z}$-invariant.  

The third  definition of the triple ratio is borrowed from (\cite{Gon2}, Section 4.2). Every point (resp. line) of $\bbRP^2$ corresponds to a line (resp. plane) in $V_3$, and we shall denote them by the same letters. The line $C$ belongs to the plane $c$ and defines a linear map $C:c\cap a\rightarrow c\cap b$: it is the graph of this map. Similarly there are linear maps $A:a\cap b\rightarrow a\cap c$ and $B:b\cap c\rightarrow b\cap a$. The composition of these three linear maps is the multiplication by the invariant $X$. 

Elaborating the third definition using a Euclidean structure, we come to the fourth definition, 
 making connection to the classical Ceva and Menelaus theorems:
$$
X:= \pm\frac{|A(a\cap b)||B(b\cap c)||C(c\cap a)|}{|A(a\cap c)||B(b\cap a)||C(c\cap b)|}, 
$$
where the distances are measured with respect to any Euclidean structure on $\mathbb R^2=\mathbb{RP}^2-\mathbb{RP}^1$ containing the triangles and the sign is positive if the triangles are inscribed on into another and minus otherwise.  The Ceva theorem claims that the lines $A(b\cap c)$, $B(c\cap a)$ and $C(a\cap b)$ are concurrent (pass through one point) if and only if $X=1$. The Menelaus theorem claims that the points $A,B$ and $C$ are concurrent (belong to a line) if and only if $X=-1$. 

In particular this definition implies the
\begin{lemma}\label{positive3}
$X$ is positive if and only if $ABC$ is inscribed into $(a\cap b)(b\cap c)(c\cap a)$.
\end{lemma}

\begin{figure}[ht]
\centerline{\epsfbox{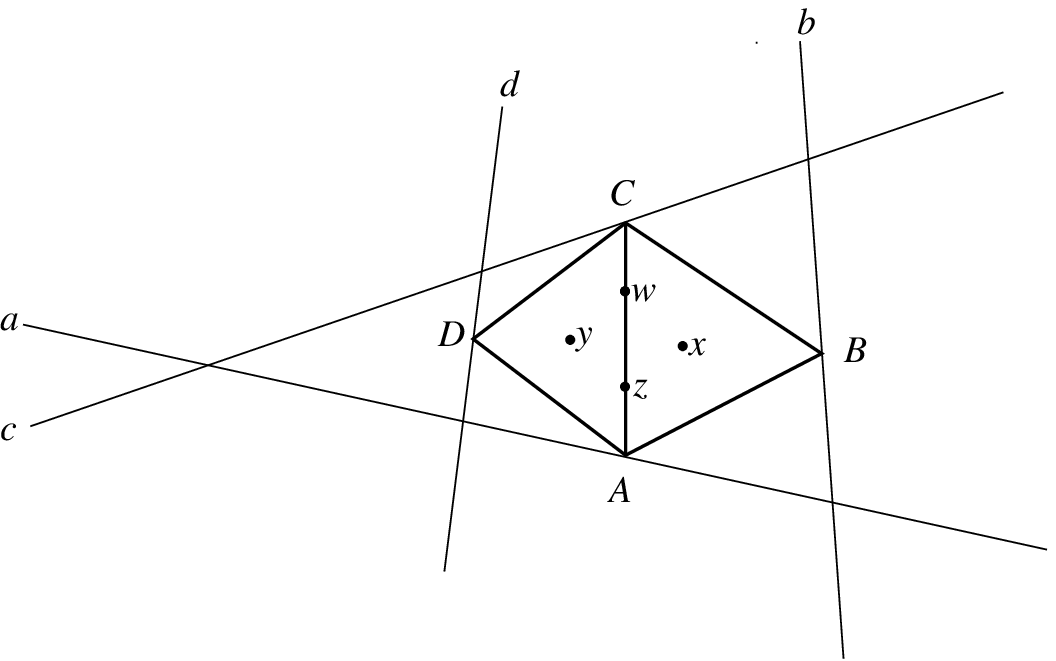}}
\caption{${\cal P}_3^4$}
\label{quadrangle}
\end{figure}

Now consider the next case:  ${\cal P}_3^4$. It is the space of pairs of quadrilaterals $abcd$ and $ABCD$, where the second is inscribed into the first, considered up to projective transformations. The space of such configurations has dimension four. Two parameters of a configuration are given by the triple ratio $X$ of the triangles $ABC$ inscribed into $abc$ and the triple ratio $Y$ of $ACD$ inscribed into $acd\{B\}$ (Figure \ref{quadrangle}). Another two parameters are given by the cross-ratios of quadruples of lines $a,AB,AC,AD$ and $c,CD,CA,CB$, denoted by $Z$ and $W$, respectively. Assign the coordinates $X,Y,Z$ and $W$ to the marked points as shown on Figure \ref{quadrangle}.

\begin{lemma} \label{positive4}
The convex quadrangle $ABCD$ is inscribed into the convex quadrangle $abcd$ if and only if  $X, Y, Z, W$ are positive. 
\end{lemma}

{\bf Proof}. Indeed, Lemma \ref{positive3} implies that $X, Y>0$.  Observe that the cross-ratio of four points on $\R {\mathbb P}^1$ is positive if and only if the cyclic order of the points is compatible with the one provided by their location on $\R {\mathbb P}^1$, i.e. with one of the two orientations of  $\R {\mathbb P}^1$. This immediately implies that $Z, W>0$. On the other hand, given $X>0$ we have, by Lemma \ref{positive3}, a unique projective equivalence class of a triangle $ABC$ inscribed into a triangle $abc$. We use the coordinates $Z, W$ to define the lines $AD$ and $CD$, and the positivity of $Z$ and $W$ guarantees that these lines intersect in a point $D$ located inside of the triangle $abc$, that is in  the same connected component   as the triangle $ABC$. Finally, $Y$ is used to define the line $d$, and positivity of $Y$ just means that it does not intersect the interior of the quadrangle $ABCD$. The lemma is proved. 

Now let us consider the space ${\cal P}_3^n$ for an arbitrary $n>2$. A point of this space is represented by a polygon $A_1\ldots A_n$ inscribed into $a_1\ldots a_n$. Cut the polygon $A_1 \ldots A_n$ into triangles. Each triangle $A_iA_jA_k$ of the triangulation is inscribed into $a_i,a_ja_k$, and we can assign their triple ratio to the marked point inside  of $A_iA_jA_k$. Moreover for any pair of adjacent triangles $A_iA_jA_k$ and $A_jA_kA_l$ forming a quadrilateral $A_iA_jA_kA_l$ inscribed into the quadrilateral $a_ia_ja_ka_l$ one can also compute a pair of cross-ratios and assign it to two points on the diagonal $A_jA_k$. The cross-ratio of the lines passing through $A_j$ is associated to the point closer to the point $A_j$.

The converse construction is straightforward. For every triangle of the triangulation one can construct a pair of triangles using the number in the center as a parameter and using Lemma \ref{positive3}. Then we assemble them together using numbers on edges as gluing parameters and Lemma \ref{positive4}. The theorem is proved.

Now let us proceed to convex projective structures on surfaces.

\vskip 3mm
\paragraph{4. A set of global coordinates on ${\cal T}_3^+(S)$.}  
Let $S$ be a Riemann surface of genus $g$ with $s$ boundary components.
Assume, that $s \geq 1$ and moreover $s\geq
3$ if $g=0$. Shrink all boundary components to points. Then the surface $S$ can be
cut into triangles with vertices at the shrunk boundary components. 
We call it an {\it ideal triangulation} of $S$. Let
us put two distinct marked points to each edge of the triangulation
and one marked point to the center of every triangle.
\begin{theorem} \label{iii}
Given an ideal triangulation $T$ of $S$, 
there exists a canonical isomorphism
$$
\varphi_{T}: {\cal T}_3(S) \stackrel{\sim}{\lra} {\R}_{>0}^{\mbox{\{marked points\}}}
$$  
\end{theorem}

{\bf Proof}. It is quite similar to the proof for ${\cal P}_3^n$. 
Let us first construct a surface starting from a collection of positive real numbers
 on the triangulation. Consider the universal cover $\widetilde{S}$ of the surface 
and lift the triangulation together with the marked points and numbers from $S$ to it. 
According to Theorem \ref{1.20.04.5}  
to any finite polygon composed of the triangles of the arising infinite 
triangulated polygon we can associate a pair of polygons in $\bbRP^2$. Let 
$U$ be the 
union of all inscribed polygons corresponding to 
such finite sub-polygons. Observe that 
it coincides with the 
intersection of all circumscribed polygons, and therefore is convex. 
The group $\pi_1(S)$ acts on 
$\widetilde {S}$ and hence on $U$ by projective transformations. The desired 
projective surface is $U/\pi_1(S)$. 

Now let us describe the inverse map, i.e. construct the 
 numbers out of a given framed convex projective structure and a given triangulation. 
Take a triangle and send it to $\bbRP^2$ using a developing map $\varphi$. 
The vertices of the triangle correspond to boundary components. For a given boundary 
component $C_i$ of $S$ 
the choice of the developing map $\varphi$ and the framing allows to 
assign a canonical flag $(A,a)$ on $\bbRP^2$  
invariant under the action of the monodromy operator around $C_i$. 
Indeed, if the boundary component $C_i$ is non-degenerate, we take the line containing the interval $\varphi(C_i)$ for $a$, and one of the endpoints on this interval for $A$. The choice between the two endpoints is given by the framing, so that the interval is oriented out from the point $A$. If the boundary component is degenerate the point $A$ is just the image of $C_i$ under $\varphi$. The line $a$ is the projectivization of the two-dimensional subspace where the monodromy operator $\mu(C_i)$ acts by a unipotent transformation.   

Assigning the flags to all three vertices of the triangle one gets a point of ${\cal P}^3_3$. The corresponding coordinate is associated to the central marked point of the original triangle. Similarly, taking two adjacent triangles of the triangulation one obtains the numbers for the marked points on edges.

These two constructions are evidently mutually inverse to each other. The theorem is proved. 

\vskip 3mm
\paragraph{5. Properties of the constructed coordinates for both ${\cal P}^3_3$ and ${\cal T}_3(S)$.}~\\

\begin{figure}[ht]
\centerline{\epsfbox{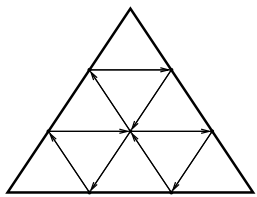}}
\caption{Poisson structure tensor.}
\label{poisson}
\end{figure}

{\bf 1.}
It turns out that the Poisson brackets between the coordinates are very simple, namely 
\begin{equation} \label{1.29.04.10}
\{X_i,X_j\}= 2\varepsilon_{ij} X_iX_j, 
\end{equation} 
where $\varepsilon_{ij}$ 
is a
skew-symmetric integral valued function. 
To define the function $\varepsilon_{ij}$ consider the graph with
vertices in marked points and oriented edges 
connecting them as shown on Figure \ref{poisson}. (We have shown edges
connecting marked points belonging to one 
triangle only. Points of other triangles are connected by arrows in
the same way. Points connected by edges without 
arrows are not taken into account.) Then 
$$
\varepsilon_{ij}=(\mbox{number of arrows from $i$ to $j$})- (\mbox{number of arrows from $j$ to $i$})
$$
The proof of this formula amounts to a long computation. However
 it is an easy exercise to verify that this bracket does not depend on the 
triangulation, and therefore can be taken just as a definition. For ${\cal P}^n_3$ it is the only natural way to define the Poisson structure known to us. 

\begin{figure}[ht]
\centerline{\epsfbox{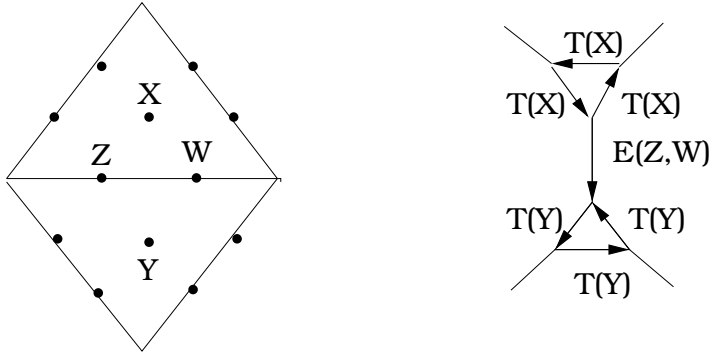}}
\caption{Construction of the monodromy group.}
\label{matrix}
\end{figure}

{\bf 2.}
Once we have positive numbers assigned to marked points, 
the construction of the corresponding projective surface is explicit. 
In particular one 
can compute the monodromy group of the corresponding projective structure
 or, in other words, the image of the map $\mu$. To describe 
the answer we use the following picture. Starting from the triangulation, 
 construct a graph embedded into the surface by drawing small edges 
transversal to each side of the triangles and inside each triangle 
connect the ends of edges pairwise by three more edges, as shown on Figure  
\ref{matrix}. Orient the edges of the triangles in counterclockwise 
direction and the other edges in the arbitrary way. Let 
$$
T(X)=\left(\begin{array}{ccc}
0&0&1\\
0&-1&-1\\
X&1+X&1
\end{array}\right)
\mbox{ and }
E(Z,W)= \left(\begin{array}{ccc}
0&0&Z^{-1}\\
0&-1&0\\
W&0&0
\end{array}\right).
$$
Assign the matrix $T(X)$ to the edges of each triangle with $X$ assigned to its center. And assign the matrix $E(Z,W)$ to the edge connecting two triangles, where $W$ (resp. $Z$) is the number to the left (resp. right) from this edge (Figure \ref{matrix}). Then for any closed path on the graph one can assign an element of $PSL_3({\mathbb R})$ by multiplying the group elements (or their inverses if the path goes along the edge against its orientation) assigned to edges the path passes along. 
The image of the fundamental group of the graph is just the desired monodromy group.


The proof of this statement is also constructive. Once we have a configuration of flags from ${\cal P}^3_3$ with triple ratio $X$, one can define a projective 
coordinate system on $\bbRP^2$. Namely, take the one where the point $b\cap c$ has coordinates $[0:1:0]^t$, the point $A$ has coordinates $[1:-1:1]^t$, $C$ -- $[1:0:0]^t$ and $B$ -- $[0:0:1]^t$ (Figure \ref{triangle}). The line $a$ has coordinates$ [1:1+X:X]$. The cyclic permutation of the flags induces the coordinate change given by the matrix $T(X)$. Similarly if we have a quadruple of flags $F_1F_2F_3F_4$ with two cross-ratios $Z$ and $W$, then the coordinate system related to the triple $F_2F_4F_1$ is obtained from the coordinate system related to the triple $F_4F_2F_3$ by the coordinate change $E(Z,W)$
.

{\bf 3.} The involution $\sigma$ acts in a very simple way, 
where $Z'=\frac{W(1+Y)}{Y(1+X)}$ and $W'=\frac{Z(1+X)}{X(1+Y)}$: 
\begin{figure}[ht]
\centerline{\epsfbox{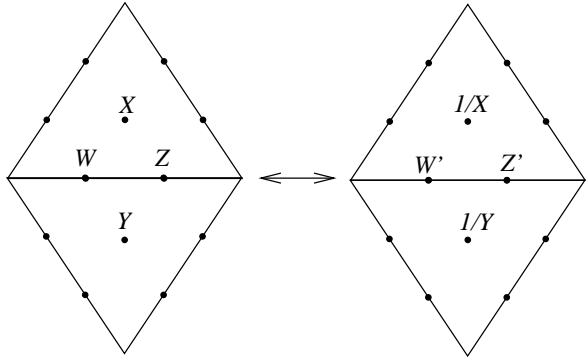}}
\caption{Involution $\sigma$.}
\label{involution}
\end{figure}


In particular a point of ${\cal T}_3(S)$ is stable under $\sigma$ if 
the two coordinates on each edge coincide, and the coordinates in the 
center of each triangle are equal to one. Taking into account that the set 
of $\sigma$-stable points is just the ordinary Teichm{\"u}ller space, one 
obtains its parametrisation. 
Actually it coincides with the one described in \cite{Teich}.

\begin{figure}[ht]
\centerline{\epsfbox{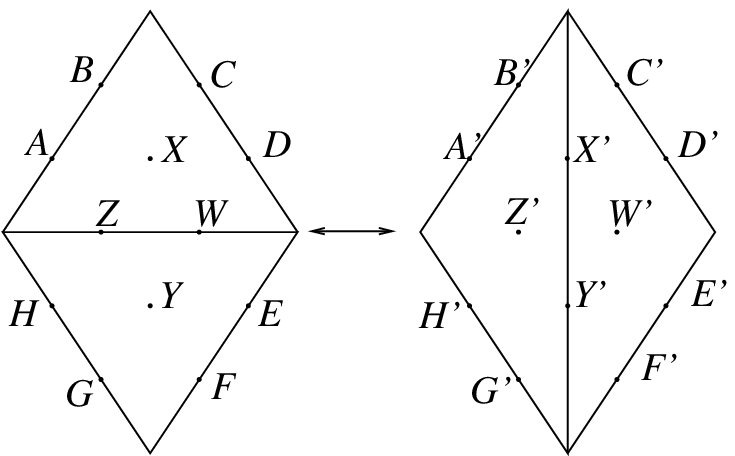}}
\caption{Flip.}
\label{flip}
\end{figure}

{\bf 4.}
Each triangulation of $S$ provides its own coordinate system and in general 
the transition from one such system to another one is given by complicated 
rational maps. However any change of triangulation may be decomposed into a 
sequence of elementary changes --- the so called {\em flips}. A flip removes 
an edge of the triangulation and inserts another one into the arising quadrilateral 
as shown on Figure \ref{flip}. This figure shows also how the numbers at the marked 
points change under the flip. Observe 
 that these formulae allow in particular to pass from one triangulation to the 
same one, but moved by a nontrivial element of the mapping class group of $S$, 
and thus give explicit formulae for the mapping class group action.


Here 
$$
A'=A(1+Z),  ~~~~D'=D\frac{W}{1+W}, ~~~~E'=E(1+W), 
~~~~H'=H\frac{Z}{1+Z},
$$
$$
B'=B\frac{1+Z+ZX+ZXW}{(1+Z)},
~~~~C'=C\frac{(1+W)XZ}{1+Z+ZX+ZXW}
$$
$$
F'=F\frac{1+W+WY+WYZ}{(1+W)}, ~~~~G'=G\frac{(1+Z)YW}{1+W+WY+WYZ},
$$
$$
X'=\frac{1+Z}{XZ(1+W)}, ~~~~Y'=\frac{1+W}{YW(1+Z)},
$$
$$
Z'=X\frac{1+W+WY+WYZ}{1+Z+ZX+ZXW},
~~~~W'=Y\frac{1+Z+ZX+ZXW}{1+W+WY+WYZ}.
$$
The formulae can be derived directly, 
or more simply as in Section 11 in \cite{Higher}. We will return to discussion of the structure of these 
formulas in Section 2.5 below.

\begin{figure}[ht]
\centerline{\epsfbox{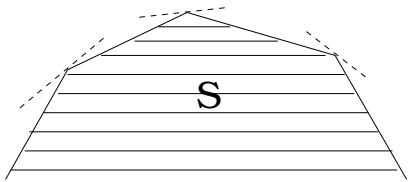}}
\caption{Framed piecewise geodesic boundary of a 
convex projective structure on S}
\label{line}
\end{figure}

\vskip 3mm
\paragraph{6. Convex projective structures on surfaces with piecewise geodesic boundary.}  
A  convex projective structure on a surface has  
{\it piecewise geodesic} boundary if  
a neighborhood of every boundary point is projectively 
isomorphic to  a neighborhood of a boundary point  of a half plane or  
 a vertex of an angle. In the latter case 
it makes sense to consider a line passing through the vertex ``outside'' of 
the surface, see the punctured lines on Figure \ref{line}.  

Let us introduce {\it framed convex  projective structure  
with piecewise geodesic boundary}.  Let $\widehat S$ be a pair consisting of 
an oriented surface $S$ with the boundary $\partial S$, and a finite 
(possibly empty) 
collection of {\it distinguished} points $\{x_1, ..., x_k\}$ on the 
boundary.   
We define the {\it punctured boundary} of $\widehat S$ by $$
\partial \widehat S:= \partial  S -\{x_1, ..., x_k\}
$$ 
So the connected components of $\partial \widehat S$ are 
either circles or arcs.

\begin{definition} \label{qwe} A framed convex  projective structure  
with piecewise geodesic boundary on $\widehat S$ is 
a convex projective structure on $S$ with the following data at the boundary:

i) Each connected component of the punctured 
boundary $\partial \widehat S$ is a geodesic interval. Moreover  
at each point $x_i$ we choose a line passing 
through $x_i$ outside of $S$, as on Figure \ref{line}.

ii) If a boundary component does not contain 
distinguished  points, we 
choose its orientation. 

The space ${\cal T}_3(\widehat S)$  parametrises framed 
convex projective structures 
with piecewise geodesic boundary on $\widehat S$. 
\end{definition}

  {\bf Remark}. Sometimes it is useful to employ a different 
definition of the framing. Namely, instead of choosing lines through 
the distinguished points $x_i$ we  choose a point 
$p_i$ inside of each of the geodesic segments bounded by the distinguished points. 
${\cal T}_3(\widehat S)$ is the moduli space
of each of these two structures on $S$. 
Indeed, the equivalence between the two definitions is seen as follows. 
Consider the convex hull of the points $p_i$. We get a surface $S'$
inscribed into $S$, with the induced convex 
projective structure. Conversely, given $S'$
we reconstruct $S$ by taking the circumscribed surface.

The convexity implies that each boundary component develops into a convex 
infinite polygon connecting two points in $\R {\Bbb P}^2$, 
stabilized by the monodromy around the component. Moreover 
at each vertex of this polygon there is a 
line segment located outside of $\widetilde S$. 

Both spaces ${\cal T}_3(S)$ and ${\cal P}^n_3$ are 
particular cases of the moduli space ${\cal T}_3(\widehat S)$.
Indeed,  ${\cal T}_3(S)$ is obtained when the set $\{p_1, ..., p_k\}$ 
is  empty, and we get ${\cal P}^n_3$ when $S$ is a disc: in this 
case $S$ is a convex polygon, serving as 
the inscribed polygon, and  the circumscribed polygon is 
given by the chosen lines.

To introduce the coordinates let us 
shrink the boundary components without marked points into punctures,  
and take a triangulation $T$ of $S$ 
with vertices at the distinguished points on the boundary and at the 
shred boundary components. We call it  {\it  ideal triangulation}. 
The interior of each  side of such a triangulation is either 
inside of $S$, or at the boundary. We assign the coordinates 
to the centers of the triangles and 
pairs of marked points on the internal sides of the triangulations. 
Namely, going to the universal cover $\widetilde S$ of $S$, we get a 
triangulation $\widetilde T$ there. 
Then, by the very definition, for every  vertex $\widetilde p$ of the triangulation 
$\widetilde T$ there is a flag given at $\widetilde p$.
We use these flags just as above.  

The space ${\cal T}_3(\widehat S)$ has the same features  1.-4. as 
${\cal T}_3( S)$. The only thing deserving a comment is 
construction of a canonical embedding 
$i: {\cal T}_2(\widehat S)\hra {\cal T}_3(\widehat
S)$. 
We define the Teichm{\"u}ller space ${\cal T}_2(\widehat S)$ 
as the space of pairs (a complete 
hyperbolic metrics on $S$, a distinguished collection of points 
$\{x_1, ..., x_k\}$ located at the absolute of $S$). A point of the absolute
of $S$ can be thought of as a  $\pi_1(S)$-orbit on the absolute of the
hyperbolic plane which is at the infinity of an end of $S$. 
Thus going to the universal cover of $S$ in the Klein model 
we get a finite collection of distinguished points at every arc of the
absolute corresponding to an end of $S$. It remains to make, in a
$\pi_1(S)$-equivariant way, finite geodesic polygons
with vertices at these points, and add the tangent 
lines to the absolute at
the 
vertices of these polygons. 
Taking the quotient by the action of $\pi_1(S)$ we get a framed convex
projective structure on $S$ with piecewise geodesic boundary.

\vskip 3mm
\paragraph{7. Comparing the moduli space of local systems and of convex 
projective structures.} Let $F$ be a set. Recall that an {\it $F$-local system} 
on a topological space $X$ is a locally trivial bundle ${\cal F}$ on $X$, 
with fibers isomorphic to 
 $F$ and 
 locally constant transition functions. If $X$ is a manifold, this is the same 
thing as a bundle with a flat connection on $X$. 
If $F$ is a group, we also assume that it acts 
on the right on ${\cal F}$, and the fibers are principal homogeneous $H$-spaces. 
If $G$ is an algebraic group, a $G$-local system on a variety $X$ 
is an algebraic variety such that for a field $K$ the set of its $K$-valued points 
is a $G(K)$-local system on $X(K)$.  
 
 Let $G := PGL_3$ and let ${\cal B}$ be the
corresponding flag variety. The set of its real points 
of was denoted by $F_3$ above. 
For a $G$-local system ${\cal L}$ on $S$ let ${\cal L}_{\cal B}:= 
{\cal L} \times_G{\cal B}$ be the associated local system of flags.  
Recall (Section 2 of \cite{Higher}) that a
 {\it framed $G$-local system} on $\widehat S$ is a pair $({\cal L}, \beta)$ 
consisting of a $G$-local system on $S$ and a flat section $ \beta$ of the
restriction of ${\cal L}_{\cal B}$ to the punctured boundary $\partial
\widehat S$. 
The corresponding moduli space is denoted by ${\cal X}_{G, \widehat S}$. 

Let us shrink all holes on $S$ without distinguished points to punctures. 
An {\it ideal triangulation} $T$ of $\widehat S$ is a triangulation 
of $S$ with vertices either at punctures or at the 
boundary components, so that 
each connected component of $\partial \widehat S$ hosts exactly one vertex of the triangulation. 
Repeating the construction 
of Sections 2.4-2.5 we get a set of coordinate systems on 
${\cal X}_{G, \widehat S}$ parametrized by ideal triangulations of
$\widehat S$. Precisely, given ideal triangulations $T$ there is an open
embedding 
$$
\varphi_T: {\mathbb G}_m^{\{\mbox{internal edges of $T$}\}} \hra {\cal X}_{G, \widehat S}
$$
 Here ${\mathbb G}_m$ stays for the multiplicative group 
understood as an algebraic group. So the set of its complex points is $\C^*$. 
The
 natural coordinate on ${\mathbb G}_m$ provides a rational coordinate function 
on the moduli space ${\cal X}_{G, \widehat S}$. Thus the map $\varphi_T$ 
provides a (rational) coordinate system on the moduli space. 
The transition functions between the two coordinate 
systems corresponding to triangulations related by a flip 
are given by the formulas written in the very end of Section 2.6. 
Since these formulas are subtraction free, the real locus of  ${\cal
 X}_{G, \widehat S}$  contains a well-defined  
subset 
$$
{\cal X}_{G,
 \widehat S}(\R_{>0}):= \varphi_T\left(\R_{>0}^{\{\mbox{internal edges of $T$}\}}\right)
\hra {\cal X}_{G, \widehat S}(\R)
$$

\begin{theorem}\label{1.24.04.1}
There is a canonical identification ${\cal X}_{PGL_3, \widehat S}(\R_{>0}) = 
{\cal T}_{3}(\widehat S)$.
\end{theorem}
In other words, ${\cal T}_{3}(\widehat S)$ is the subset 
of the real locus of the moduli space ${\cal X}_{PGL_3, \widehat
  S}$ 
consisting of the points which have positive real coordinates in the
coordinate system corresponding to one (and hence any) ideal triangulation of 
 $\widehat S$. 

{\bf Proof}. We will assume that ${\cal T}_{3}(\widehat S)$ parametrises the
structures 
described in the Remark after Definition \ref{qwe}. (If we adopt  
Definition \ref{qwe},   
the distinguished points will 
play different roles in the definition of these moduli spaces). 
Let us define a canonical embedding ${\cal T}_{3}(\widehat S) \to 
{\cal X}_{G, \widehat
  S}(\R_{>0})$. The local system ${\cal L}_p$ corresponding to a  point $p \in {\cal
  T}_{3}(\widehat S)$ 
is the one corresponding to the representation $\mu(p)$. Intrinsically, 
it is
the local system of projective frames on $S$. Let us define a
framing  on it. Let $C_i$ be a boundary component of $S$. If it does not contain 
the distinguished points, the corresponding component of the framing is provided
by the flag $(a, A)$ defined in the proof of Theorem \ref{iii}. If $C_i$ 
contains distinguished points, we use the flag $(p, L)$ where $p$ is the
chosen point
on the boundary geodesic interval $L$. The restriction of the canonical coordinates 
on the moduli space ${\cal X}_{PGL_3, \widehat
  S}$ to the space of convex projective structures gives,  by the very 
definition, the coordinates on the latter space defined above. 
Now the proof follows immediately from the results of Sections 2.3-2.5. 
The theorem is proved.  

\vskip 3mm \paragraph{8. 
Laurent positivity properties of the monodromy representations.}   
Consider the universal $PGL_3$-local system  on
the space $S \times {\cal X}_{PGL_3, S}$. 
Its fiber over $S \times p$ is isomorphic to 
the local system corresponding to the
point $p$ of the moduli space ${\cal X}_{PGL_3, S}$.
Our
next goal is to study its  monodromy 
representation. 

Let ${\mathbb F}_S$ be the field of rational functions on the moduli
space ${\cal X}_{PGL_3, S}$. The
monodromy of  the universal local system around a loop on $S$ is a conjugacy class in
$PGL_3({\mathbb F}_S)$. 

Recall that, given an ideal triangulation $T$ of $S$, we defined 
canonical coordinates $\{X^{T}_i\}$ on ${\cal X}_{PSL_3, S}$
corresponding to $T$. Given a coordinate system $\{X_i\}$, 
a {\it positive rational function in $\{X_i\}$} 
is a function which can be presented as a ratio 
of two polynomials in $\{X_i\}$ with positive integral coefficients. 
For instance $1-x+x^2$ is a positive rational function since it is equal to 
$(1+x^3)/(1+x)$.  
A matrix 
is {\it totally positive} if all its minors are non-zero positive
rational functions  in $\{X^{T}_i\}$. 
Similarly we define upper/lower triangular totally positive 
matrices: any minor, which is not identically zero for generic
matrix in question, is a positive
rational functions.  

\begin{theorem} \label{2.1.04.0} Let $T$ be an ideal
triangulation of $S$. 
Then  monodromy of the universal $PGL_3$-local system on $S \times {\cal
  X}_{PGL_3, S}$ around 
any non-boundary loop on $S$ is conjugate in $PGL_3({\mathbb
  F}_S)$  to a totally
positive matrix in  the
canonical coordinates assigned to $T$.

The monodromy around a boundary component is conjugated to an
upper/lower 
triangular totally
positive matrix. 
\end{theorem}

{\bf Proof}. Recall the trivalent graph ${\Gamma}'_T$ on Figure 4,
defined starting from an ideal
triangulation $T$ of $S$. Let us call its edges forming the little
triangles by $t$-edges, and the edges dual to the edges of the
triangulation $T$  by  $e$-edges. The matrices $T(*)$ and $E(*)$ assigned in
Section 5.2 to the $t$- and $e$-edges of this graph provide an
explicit construction of the universal local system. Indeed, since 
$T(X)^3$ is the identity in $PGL_3$, we get a universal local system
on the dual graph to $T$. This graph is homotopy equivalent to $S$. 

Given a loop on $S$, we shrink it to a loop on
${\Gamma}'_T$. We may assume that this loop contains no 
consecutive $t$-edges: indeed, a composition of two $t$-edges is a
$t$-edge. Thus we may choose an initial vertex on the loop so that the
edges have the pattern $et ... etet$. Therefore the monodromy is computed as a
product of matrices of type  $ET$ or $ET^{-1}$. Observe that   
\begin{equation} \label{3.18.04.3}
E(Z,W) T(X)=\left(\begin{array}{ccc}
Z^{-1}X &Z^{-1}(1+X) &Z^{-1} \\
0 & 1 & 1\\
0 & 0 & W^{-1}
\end{array}\right)\quad \mbox{and}
\end{equation}
\begin{equation} \label{3.18.04.4}
E(Z,W)T(X)^{-1}= \left(\begin{array}{ccc}
Z^{-1}&0&0\\
1&1 &0\\
W& W(1+X^{-1})&WX^{-1} 
\end{array}\right)\quad \mbox{}.
\end{equation}
These matrices 
are upper/lower 
triangular totally
positive 
integral Laurent matrices in the coordinates $X, Y, Z, W$. 
A loop on $S$ contains matrices of just one kind if and only if it is the loop
around a boundary component. Therefore the monodromy around any non-boundary
loop is obtained as a product of both lower and upper
triangular matrices. It remains to use the following fact: 
if $M_i$ is either lower or upper triangular totally positive matrix,
and there are both upper and lower triangular matrices among $M_1,
..., M_k$, then the product $M_1 ... M_k$ is a totally positive 
matrix (cf. \cite{GKr}). The theorem is
proved. 

We say that the monodromy of a $PGL_3(\R)$-local system is {\it regular
  hyperbolic} if the monodromy around a non-boundary loop 
has distinct real  eigenvalues, and the monodromy around a boundary
  component is conjugated to a real  totally positive upper
  triangular matrix. 

\begin{corollary} \label{2.1.04.1} 
The  monodromy of a convex projective structure with geodesic boundary
on $S$ is faithful and 
regular hyperbolic. 
\end{corollary}

{\bf Proof}. Follows immediately from Theorem \ref{2.1.04.0} and the
Gantmacher-Krein theorem \cite{GKr} (see Chapter 2, Theorem 6 there) 
claiming that the eigenvalues of a totally positive matrix are distinct real numbers. The corollary is proved 
\footnote{As was pointed out by the Referee, an alternative proof of 
Corollary \ref{2.1.04.1} follows from the results of \cite{Goldman} for  surfaces without boundary by doubling along the geodesic boundary.}

{\bf Remark}. Both Theorem \ref{2.1.04.0} and
Corollary \ref{2.1.04.1} 
can be generalized to the case when $G$ is an  arbitrary split reductive group
 with trivial center, see
\cite{Higher}. When the boundary of $S$ is empty and $G= PSL_m(\R)$, 
a statement similar to Corollary \ref{2.1.04.1} 
was proved in \cite{Lab2}. 

 A Laurent polynomial in $\{X^{T}_i\}$
is {\it positive integral} if its coefficients are positive integers. 
A rational function on the moduli space ${\cal X}_{PGL_3, S}$ 
is called a {\it good positive Laurent polynomial} if it is 
a positive 
integral Laurent polynomial in every
canonical coordinate system on  ${\cal X}_{PGL_3, S}$.

\begin{corollary} \label{2.1.04.2} 
The  trace of the $n$-th power of the 
monodromy of the universal $PGL_3$-local system 
on $S$ around any loop on $S$ is a good positive Laurent polynomial 
on  ${\cal X}_{PGL_3, S}$. 
\end{corollary}

{\bf Proof}. A product of matrices with 
positive integral Laurent coefficients is again a matrix of this type. 
So the corollary follows from the explicit construction of the universal 
$PGL_3$-local system given in the proof of Theorem \ref{2.1.04.0}, 
and formulas (\ref{3.18.04.3})-(\ref{3.18.04.4}).

A good positive Laurent polynomial on ${\cal X}_{PGL_3, S}$ 
is {\it indecomposable} if it can not be presented as a sum 
of two non-zero good positive Laurent polynomials on 
${\cal X}_{PGL_3, S}$. 

\begin{conjecture} \label{2.1.04.5} 
The  trace of the $n$-th power of the 
monodromy of the universal $PGL_3$-local system 
on $S$ around any loop on $S$ is indecomposable.  
\end{conjecture}

\section{The universal  $PSL_3(\R)$-Teichm{\"u}ller space } 

In this Section, using the positive configuration of flags, 
 we introduce the universal  
Teichm{\"u}ller space ${\cal T}_3$ for $PSL_3(\R)$. 
We show that the 
Thompson group acts by its automorphisms, 
preserving a natural Poisson structure on it. We show that ${\cal T}_3$ 
is closely related to the space of convex curves in $\R {\mathbb P}^2$. 

We leave to the reader to formalize a definition of a finite cyclic set. 
As a hint we want to mention that a subset of an oriented circle inherits a cyclic structure. 

\vskip 3mm \paragraph{1. The universal  higher Teichm{\"u}ller space.}  
\begin{definition} i) A set $C$ is {\em cyclically ordered} if  
  any finite
  subset of $C$ is cyclically ordered, and any inclusion of finite subsets
  preserves the cyclic order. 

ii) A map $\beta$ from a cyclically ordered set $C$  to the flag variety
  $F_3$ is {\em positive} if it maps every cyclically ordered quadruple 
  $(a,b,c,d)$ to a positive 
  quadruple of flags $(\beta(a), \beta(b), \beta(c), \beta(d))$.
\end{definition}

A positive map is necessarily injective. So, abusing
terminology, we define {\it positive subset} $C \subset F_3$ as the image of a
positive map, keeping in mind a cyclic structure on $C$. 

The set of positive $n$-tuples of flags in $F_3$ is
invariant under the operation of reversing the cyclic order to the
opposite one. So it has not only cyclic, but also dihedral symmetry.

Consider a cyclic
order on ${\mathbb P}^1(\Q)$ provided by the 
embedding ${\mathbb P}^1(\Q) \hra {\mathbb P}^1(\R)$.

\begin{definition} \label{1.29.04.3} The universal  higher Teichm{\"u}ller space ${\cal T}_3$ 
     is the quotient of the space of positive maps
 $\beta: {\mathbb P}^1(\Q) \lra F_3$
by the action of the group $PSL_3(\R)$. 
\end{definition}

{\it Canonical coordinates on ${\cal
  T} _3$}. 
Let us choose three distinct points on ${\mathbb P}^1(\Q)$,
  called $0, 1, \infty$. Recall the Farey triangulation, understood 
as a triangulation of the hyperbolic disc with a
  distinguished oriented edge.  
Then we have  canonical identifications 
\begin{equation} \label{1.29.04.6} 
{\mathbb P}^1(\Q) = \Q \cup \infty = \{\mbox{vertices of the Farey
  triangulation}\}
\end{equation} 
The distinguished oriented edge 
goes from $0$ to $\infty$. Consider the infinite set
$$
I_3 := \{\mbox{pairs of points on each edge of the Farey
  triangulation}\} \cup
$$
$$
\{\mbox{(centers of the) triangles of the Farey triangulation}\} 
$$ 
A point of ${\cal T}_3$ gives rise to  
a positive map 
$$
\{\mbox{vertices of the Farey triangulation}\} \lra F_3
$$
considered up to the action of $PGL_3(\R)$. We assign to
every triple of flags  at the vertices of a Farey triangle
their triple ratio, and to every quadruple of flags at the vertices of a Farey 
quadrilateral the related two cross-ratios,  pictured  at the diagonal. 
We get a canonical map
\begin{equation} \label{1.29.04.7} 
\varphi_3: {\cal T}_3 \stackrel{}{\lra} \R^{I_3}_{>0}
\end{equation}

\begin{theorem} \label{1.29.04.3} The map (\ref{1.29.04.7}) is an isomorphism. 
\end{theorem}

{\bf Proof}. It is completely similar to the one of Theorem
\ref{1.20.04.5}, and thus is omitted. 

The right hand side in (\ref{1.29.04.7}) has a natural quadratic
Poisson structure given by the formula (\ref{1.29.04.10}). Using the isomorphism
$\varphi_3$, we transform  it to a Poisson structure on $\mathcal T_3$. 

 Here is another way to produce points of  $\mathcal T_3$. 
\vskip 3mm \paragraph
{2. Convex curves on $\R{\mathbb P}^2$ and positive curves in
  the flag variety $F_3$.} 
Let $K$ be a continuous closed convex curve in $\R{\mathbb P}^2$. Then for every
point $p\in K$ there is a nonempty set of lines intersecting $K$ at $p$ and such that $K$ is on one side of it. We call them {\it osculating lines at $p$}. A {\it regular convex} curve in $\R{\mathbb P}^2$ is a continuous convex curve which has exactly one osculating line at each point. Assigning to a point $p$ of a  regular convex curve $K$ the unique osculating line at $p$ we get the {\it osculating curve} $\widetilde K$ in the flag variety $F_3$. It is continuous. Moreover,  thanks  to the results of Section 2.3, it is positive, i.e. the map $K \to \widetilde K \subset F_3$
  is positive. 
\begin{figure}[ht]
\centerline{\epsfbox{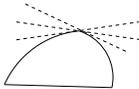}}
\caption{Osculating lines to a convex curve.}
\label{Osculating}
\end{figure}
It follows that $K$ is a $C^1$--smooth convex curve in
  $\R{\mathbb P}^2$,
and  the osculating curve $\widetilde K$ 
   is given by the 
flags $(p, T_pK)$ where $p$ runs through $K$. It turns out that this
  construction gives all positive continuous curves in the flag
  variety $F_3$. Recall the double bundle 
$$
\R {\mathbb P}^2 \stackrel{p}{\longleftarrow} F_3 \stackrel{\widehat p}{\lra} 
\R\widehat {\mathbb P}^2
$$ 
where $\R\widehat {\mathbb
 P}^2$ is the set of lines in $\R{\mathbb  P}^2$.   A subset of $\R{\mathbb P}^2$ 
is {\it strictly convex} if no line contains
its three distinct points. (warning: a strictly convex curve is a 
strictly convex subset, but 
a convex domain is never a strictly convex subset). 

\begin{proposition} \label{2.24.04.11} i) If 
$C$  a positive subset of $F_3$, then $p(C)$ and $\widehat p(C)$
are convex. 

ii) If $C$ is a 
continuous positive curve in $F_3$, then $p(C)$ is a regular convex curve. This gives rise to a
bijective correspondence 
\begin{equation} \label{2.25.04.1}
\mbox{ continuous positive curves in $F_3$} <=> \mbox{regular
convex curves in $\R{\mathbb P}^2$}. 
\end{equation}
\end{proposition}

{\bf Proof}.  i) Follows immediately from the results of Section 2.3. 

ii) $p(C)$ is continuous and, by i), convex. It comes equipped with a
continuous family of osculating lines. It follows from this that
$p(C)$ is regular. This gives the arrow $=>$ in (\ref{2.25.04.1}). 
The osculating curve provide the opposite arrow. Since a convex curve
can not have two different continuous families of osculating lines, we
got a bijection. 
The proposition is proved.

So identifying $P^1(\Q)$ with a subset
  of  a regular convex curve $K$ respecting the
  cyclic orders  
we produce a point of  ${\cal
  T}_3$. 

{\it The Thompson group ${\mathbb T}$}. It is the group of all piecewise $PSL_2(\Z)$-projective automorphisms of
${\mathbb P}^1(\Q)$. By definition,  for  every
$g \in {\mathbb T}$ there exists a decomposition of ${\mathbb P}^1(\Q)$
into a union of finite number of segments, which may intersect only at the
ends,  so that the restriction of $g$ to each segment is
given by an element of $PSL_2(\Z)$.  The Thompson group  acts
on  ${\cal T}_3 $ in an obvious way. 
Here is another way to look at it. The Thompson group contains the 
following elements, called flips: Given an edge $E$
of a triangulation $T$ with a distinguished
oriented edge, we do a flip at 
 $E$ (as on the Figure \ref{flip}), obtaining a new triangulation $T'$ with a distinguished
oriented edge. Observe that the sets of ends of the trivalent trees dual to
the triangulations $T$ and $T'$ are
identified. On the other hand, there exists unique isomorphism of the
plane trees $T$ and $T'$ which identifies  their distinguished oriented
edges. It provides a map of the ends of these trees, and hence an
automorphism of $P^1(\Q)$, which is easily seen to be piecewise
linear. The Thompson group is generated by 
flips (\cite{I}). So the formulas in the end of Section 2.5 allow to write 
the action of the Thompson group explicitly in our coordinates. 

\begin{proposition} The action of the Thompson group  preserves the
  Poisson bracket 
on  ${\cal T}_3$. 
\end{proposition} 

{\bf Proof}.  
Each flip preserves the Poisson bracket. 
The proposition follows. 

\vspace{4mm}
Let ${\cal H}$ be
the hyperbolic plane. For a torsion free subgroup $\Delta \subset 
PSL_2(\Z)$, set $S_{\Delta}:= {\cal H}/\Delta$.

\begin{proposition} \label{2.23.04.2}
The Teichm\"uller space ${\cal T}_3(S_{\Delta})$ 
is embedded into the universal Teichm\"uller space ${\cal T}_3$ as the subspace of
$\Delta$-invariants:
\begin{equation} \label{2.23.04.1}
{\cal T}_3(S_{\Delta}) = ({\cal T}_3)^{\Delta}
\end{equation} 
This isomorphism respects the Poisson brackets. 
\end{proposition} 

{\bf Demonstration}. The surface $S_{\Delta}$
has a natural triangulation $T_{\Delta}$, the image of the Farey
triangulation under the projection $\pi_{\Delta}: {\cal H} \to
S_{\Delta}$. So according to Theorem \ref{1.20.04.5} 
the left hand side in (\ref{2.23.04.1}) is identified with the
 positive
valued functions on $I_3/\Delta$, that is with $\Delta$-invariant positive
valued functions on $I_3$. It remains to use Theorem \ref{1.29.04.3}. 
The claim about the Poisson structures follows from the very
definitions. The proposition is proved. 

{\bf Remark 1.} The mapping class group of $S_{\Delta}$ is not embedded
into the Thompson group, unless we want to replace the latter by a bigger
group. Indeed, a flip at an edge $E$ on the surface 
$S_{\Delta}$ should correspond to the infinite composition of flips at the
edges of the Farey triangulations projected onto $E$. 
Nevertheless the Thompson group plays the role of the mapping
class group for the universal Teichm\"uller space. 

{\bf Remark 2.} The analogy between the mapping class groups of surfaces 
and the Thompson group 
can be made precise as follows. In Chapter 2 of \cite{Cluster} we defined 
the notion of a cluster ${\cal X}$-space (we 
recall it in the Section 4.3),  
and the 
{\it mapping class group of a cluster ${\cal X}$-space}, acting by automorphisms 
of the cluster ${\cal X}$-space. 
The classical and the universal 
Teichm\"uller spaces can be obtained as 
the spaces of $\R_{>0}$-points of certain cluster ${\cal X}$-spaces,  which were 
defined in \cite{Higher}. The corresponding mapping class groups are the classical mapping class groups, and the Thompson group, see 
Sections 2.9 and 2.14 in \cite{Cluster} for more details. 

{\bf Remark 3.} {\it The universal Teichm\"uller space and the $SL_3$ Gelfand-Dikii
  Poisson brackets}. 
The functional space of osculating curves to 
smooth 
curves has a natural Poisson structure, the $SL_3$
Gelfand-Dikii bracket. One can show that its restriction to the 
subspace provided by the regular convex curves is compatible with our Poisson
bracket.

\section{The quantum $PGL_3$-Teichm\"uller spaces }

In this Section we introduce the quantum 
universal $PGL_3$-Teichm\"uller space. We proved 
that the Thompson group acts by its automorphisms. 
In particular this immediately gives a definition of the quantum 
higher Teichm\"uller space of a punctured surface $S$, 
plus the fact that the mapping class group of $S$ acts by its automorphisms. 

\vskip 3mm
{\bf 1. The quantum torus related to the set $I_3$}. 
Let $T^q_3$ be a non-commutative algebra, (a quantum torus),
generated by the elements $X_i$, where $i \in I_3$, subject to the
relations
\begin{equation} \label{2.23.04.10}
q^{-\varepsilon_{ij}}X_iX_j = q^{-\varepsilon_{ji}}X_jX_i, \quad i,j
\in I_3
\end{equation} 
Here $q$ is a formal variable, so it is an algebra over $\Z[q,
q^{-1}]$. The algebra $T^q_3$ is equipped with an antiautomorphism
$\ast$, acting on the generators by 
\begin{equation} \label{ast}
\ast(X_i) = X_i, \quad \ast(q) =
q^{-1}
\end{equation}
 Let us denote by ${\rm Frac}(T^q_3)$ its non-commutative
fraction field. 
\vskip 3mm 
\paragraph{2.  The Thompson group action: the quasiclassical limit.}
Recall that the Thompson group ${\Bbb T}$ is generated by flips at
edges of the Farey triangulation. 
Given such an edge $E$, let us define an automorphism
$$
\varphi_E: {\rm Frac}(T_3) \lra {\rm Frac}(T_3)\qquad \mbox{where $T_3:= T^1_3$} 
$$
 Consider the $4$-gon of the Farey triangulation
 obtained by gluing the two triangles sharing the edge $E$ --- see
Figure  6, where  $E$ is the horizontal  edge. Let $S_E \subset I_3$
be the $12$-element set formed by the pairs of marked points on the
four sides and two diagonals of the $4$-gon. On Figure  6 these are the points
labeled by $A, B, ..., Z, W$. We set
$$
\varphi_E(X_i) := \left\{ \begin{array}{ll}X_i' & i \in S_E\\
X_i & i \not \in S_E
\end{array}\right.
$$
where $X_i'$ is computed using the formulas after Figure  6, so for
instance if $X_i =A$, then $X_i' =A'$ and so on.

\begin{proposition} \label{2.24.04.1}
The automorphisms $\varphi_E$ give rise to an action of the Thompson group
${\Bbb T}$ by
automorphisms of the field ${\rm Frac}(T_3)$ preserving the Poisson
bracket. 
\end{proposition}

\begin{figure}[ht]
\centerline{\epsfbox{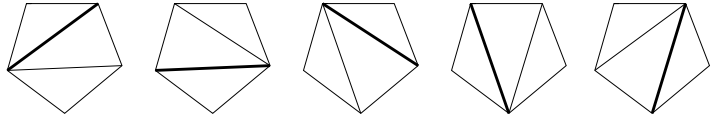}}
\caption{Pentagon relation.}
\label{Pentagonr}
\end{figure}

{\bf Proof}. It is known that the only relations between the flips
$f_E$ are the following ones:
$$
{\rm i}) f_E^2 = {\rm Id}, \qquad {\rm ii}) f_{E_1} f_{E_2} = f_{E_2} f_{E_1}
 \mbox{ if $E_1$ and $E_2$ are disjoint}, 
\qquad 
{\rm iii})(f_{E_5} f_{E_4} f_{E_3} f_{E_2} f_{E_1})^2 = {\rm Id}
$$
where the sequence of flips $f_{E_1}, ..., f_{E_5}$ is shown on Figure  \ref{Pentagonr}. 
The first two relations are obviously valid. The  pentagon relation
iii)
  is clear from the geometric origin of the
formula. The proposition is proved.

To introduce the action of the Thompson group on the non-commutative field 
${\rm Frac}(T^q_3)$, let us  recall some
facts about the quantum cluster ensembles from \cite{Cluster}. \

\vskip 3mm
\paragraph{3. The cluster ${\cal X}$-space and its quantization.} 
The cluster ${\cal X}$-space is defined by using the same set-up 
as the cluster algebras \cite{FZI}. 
We start from the quasiclassical case. Consider the following data 
$\mathbf I = (I, \varepsilon)$, 
called a seed:

i) A set $I$, possibly infinite. 

ii) A skew-symmetric function $\varepsilon = \varepsilon_{ij}: I \times I \to \Z$.

We assume that for every $i \in I$, the set of $j$'s such that
$\varepsilon_{ij} \not = 0$ is finite. 

A {\it mutation in the direction $k \in I$} 
changes the seed $\mathbf I$ to a new one 
$\mathbf I'= (I, \varepsilon')$ where the function
 $\varepsilon'$ is given by the following somewhat mysterious 
formula, which 
appeared in \cite{FZI} in the definition of
cluster algebra: 
\begin{equation}\label{5.18.04.2}
\varepsilon_{ij}' =  \left\{ \begin{array}{ll} - \varepsilon_{ij} & \mbox{ if $k
      \in \{i,j\}$}
\\\varepsilon_{ij}+ \varepsilon_{ik}{\rm max}\{0, {\rm sgn}(\varepsilon_{ik})
\varepsilon_{kj}\}& \mbox{ if $k \not 
      \in \{i,j\}$}
\end{array}\right.
\end{equation}
A seed may have an automorphism given by a bijection 
$I \to I$ preserving the $\varepsilon$-function. 
A {\it cluster transformation} is a composition of 
mutations and automorphisms of seeds.

\vskip 3mm
\noindent
We assign to a seed $\mathbf I$  a Poisson 
 algebra ${T}_{\mathbf I}:= \Q[X_i, X_i^{-1}]$, $i \in I$ with 
the Poisson bracket 
$$
\{X_i, X_j\} = \varepsilon_{ij}X_iX_j
$$
Observe that the algebra ${T}_{\mathbf I}$ is the algebra of regular functions 
on an algebraic tori ${\cal X}_{\mathbf I}$: for any field $K$, the set of its $K$-valued points is $(K^*)^I$.    

Let us take another seed $\mathbf I'= (I, \varepsilon')$ 
with the same set $I$, and  
consider a homomorphism $\mu_k: {\rm Frac}({T}_{\mathbf I'}) \lra {\rm Frac}({T}_{{\mathbf I}})$ defined on the generators by 
\begin{equation}\label{5.18.04.2as}
\mu_k(X'_i) := \left\{ \begin{array}{lll} X_i (1+X_k)^{-\varepsilon_{ik}}
 & \mbox{ if $k\not =i$, $\varepsilon_{ik} \leq 0$ }\\
X_i 
(1+X^{-1}_k)^{-\varepsilon_{ik}}& \mbox{ if $\varepsilon_{ik} > 0$ }\\
X^{-1}_k & \mbox{ if $k=i$}\end{array}\right.
\end{equation}

The following result could serve as a motivation for the 
formula  (\ref{5.18.04.2})
\begin{lemma} \label{5.18.04.1}
The map $\mu_k$ preserves the Poisson bracket if and only if the functions 
$\varepsilon'_{ij}$ and $\varepsilon_{ij}$ are related by the 
formula  (\ref{5.18.04.2})
\end{lemma}

{\bf Proof}. Straightforward. 

\vskip 3mm 
The cluster ${\cal X}$-space is a geometric object 
describing  the birational 
transformations (\ref{5.18.04.2as}). It can be understood as a 
space ${\cal X}_{|\mathbf I|}$ 
obtained by 
gluing the algebraic tori ${\cal X}_{\mathbf J}$, corresponding to all 
seeds $\mathbf J$ obtained from the initial one $\mathbf I$ by cluster 
transformations, 
according to these birational 
transformations. 

\vskip 3mm 
Now let us turn to  the quantum case. Given a seed, consider 
the 
non-commutative fractions field ${\rm Frac}({T}^q_{\mathbf I})$ 
of the algebra generated by the 
elements $X_i$, $i \in I$, subject to the relations
(\ref{2.23.04.10}). It is a non-commutative $\ast$-algebra, with an automorphism $\ast$ given by  (\ref{ast}).

It was proved in Section 3 of \cite{Cluster} that there is a $\ast$-algebra
homomorphism
$$
\mu_k: {\rm Frac}({T}^q_{\mathbf I'}) \lra {\rm Frac}({T}^q_{\mathbf I})
$$
acting on the generators by
$$
\mu_k(X_i) = \left\{ \begin{array}{lll} X_i G_{|\varepsilon_{ik}|}(q;
    X_k) & \mbox{ if $k\not =i$, $\varepsilon_{ik} \leq 0$ }\\
X_i 
G_{|\varepsilon_{ik}|}(q; X_k^{-1})^{-1}& \mbox{ if $\varepsilon_{ik} > 0$ }\\
X^{-1}_k & \mbox{ if $k=i$}\end{array}\right.
$$
where
$$
G_a(q; X) := \left\{ \begin{array}{ll}\prod_{i=1}^a(1+q^{2i-1}X) & a>0\\
1 & a=0
\end{array}\right.
$$ 
\vskip 3mm
\paragraph{4. The main result.} Our goal is the following
theorem. 
\begin{theorem} \label{2.24.04.2}
The Thompson group
${\Bbb T}$ acts by $\ast$-automorphisms of the non-commutative 
field ${\rm Frac}(T^q_3)$, so that specializing $q=1$ we recover the
action from Proposition \ref{2.24.04.1}.  
\end{theorem}

{\bf Proof}. To get the explicit
formulae for the action of flips, 
 we will apply the above construction  to the situation when $I = I_3$, and the function
$\varepsilon_{ij}$ was defined in Section 2.5. 

We define a flip
$$
\varphi_E: {\rm Frac}({T}^q_3) \lra {\rm Frac}({T}^q_3)
$$
as the composition of the mutations in the directions of $Z,
W, X, Y$ on Figure  6. An easy computation, presented in the subsection 4.6 
in the quantum form,  shows that in the case $q=1$ this leads to the formulas
after Figure  6. In the general case one needs to check the relations for
the flips. The first two of them are obvious. The quantum pentagon relation
can be 
 checked by a tedious computation, 
using the explicit formulas for the quantum flip given in the subsection 4.6. 
Observe that $f_{E_5} f_{E_4} f_{E_3} f_{E_2} f_{E_1}$ differs from the
identity only by an involution of the seven variables assigned to
the marked points inside of the pentagon. 
The theorem is proved.

\vskip 3mm 
The universal $PGL_3$-Teichm\"uller space 
story from this point of view looks as follows: 

\begin{theorem} The universal 
Teichm\"uller space ${\cal T}_3$ is the set of the positive real
points of the cluster ${\cal X}$-space related to the
seed 
$$
\{ \mbox{\rm the set $I_3$, the cluster function $\varepsilon_{ij}$ defined in
Section 2.5} \}
$$  
The Thompson group is a subgroup of the mapping class group of the corresponding 
cluster ensemble. 
\end{theorem}

{\bf Remark}. {\it The quantum universal higher 
Teichm\"uller space and the $W_3$-algebra}. 
We suggest that the action of the Thompson group by birational
$\ast$-automorphisms of the quantum universal Teichm\"uller space
${\cal T}_3^q$ should be considered as 
 an incarnation of the $W$-algebra corresponding to
$SL_3$. (For a down-to-earth definition of $W_3$-algebras see 
\cite{Belavin}; for an algebraic-geometric 
discussion of $W$-algebras see \cite{BD}).

\vskip 3mm
\paragraph{5. The quantum $PGL_3$-Teichm\"uller space for a punctured surface.} 
Let us present a punctured surface $S$ as $S = S_{\Delta}:= {\cal
  H}/\Delta$, where $\Delta$ is a subgroup of $PSL_2(\Z)$, as in Section 3.2. 

\begin{definition} \label{qSL_3} The quantum higher 
Teichm\"uller space ${\cal X}^q_3(S_{\Delta})$ for a punctured surface 
$S_{\Delta}$ is given by 
\begin{equation} \label{2.25.04.3}
{\cal X}^q_3(S_{\Delta}):= {\rm Frac}({T}^q_3)^{\Delta}
\end{equation}
\end{definition}

Theorem \ref{2.24.04.2} immediately implies 

\begin{corollary} 
The mapping class group of $S$ acts by positive birational
$\ast$-automorphisms 
of the $\ast$-algebra (\ref{2.25.04.3}).
 \end{corollary}

\vskip 3mm
\paragraph{6. Appendix: Formulas for a quantum flip.}  Performing mutations at $Z$
and $W$ we get 
$$
B_1 = B, \quad C_1 = C, \quad F_1 = F, \quad G_1 = G, \quad Z_1 = Z^{-1},
\quad W_1 = W^{-1}
$$
$$
A_1 = A(1+qZ), \quad D_1 = D(1+qW^{-1})^{-1}, \quad H_1 = H(1+qZ^{-1})^{-1},
\quad E_1 = E(1+qW),
$$
$$
X_1 = X(1+qZ^{-1})^{-1}(1+qW), \quad Y_1 = Y(1+qZ)(1+qW^{-1})^{-1}
$$
\begin{figure}[ht]
\centerline{\epsfbox{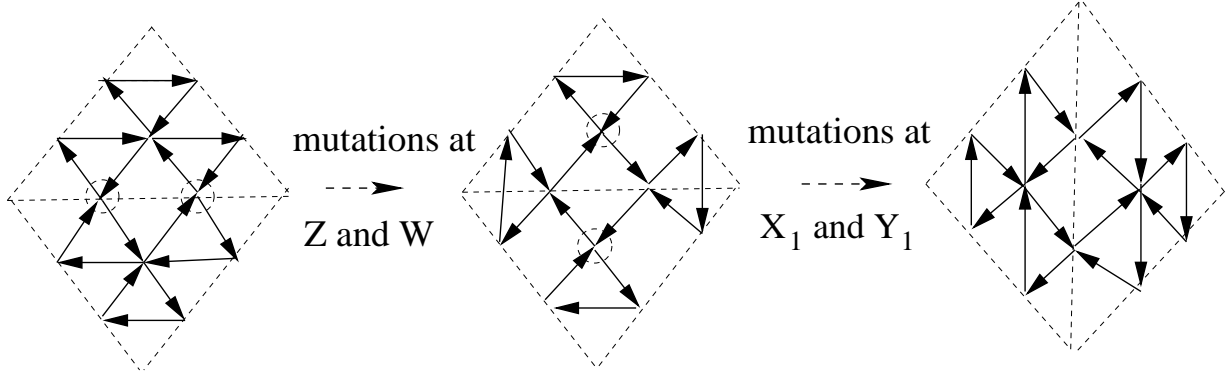}}
\caption{Mutations of the $\varepsilon$-function }
\label{Pentagon}
\end{figure}
Mutations of the $\varepsilon$-function  are illustrated on Figure \ref{Pentagon}, 
where the centers of mutations are shown by little circles.

Performing mutations at $X_1$ and $Y_1$ we arrive at the final
formulas for the quantum flip:
$$
A' = A(1+qZ), \quad D' = D(1+qW^{-1})^{-1}, \quad E' = E(1+qW), 
\quad H' = H(1+qZ^{-1})^{-1}, 
$$
$$
X'= (1+qW)^{-1}(1+qZ^{-1})X^{-1}, \quad Y'= (1+qW^{-1})(1+qZ)^{-1}Y^{-1}, 
$$
$$
B' =  B (1+qZ)^{-1}(1+qZ + q^2ZX + q^3ZXW)
$$
$$
F' =   F(1+qW)^{-1}(1+qW + q^2WY + q^3WYZ)
$$
$$
C' = CZX(1+qW)\Bigl(1+q^{-1}Z + q^{-2}XZ + q^{-3}WXZ\Bigr)^{-1} 
$$
$$
G' = GWY(1+qZ)(1 +q^{-1}W+q^{-2}YW + q^{-3}ZYW)^{-1}
$$
$$
Z'= 
X(1+qW) 
\Bigl(1+q^{-1}Z + q^{-2}XZ + q^{-3}WXZ \Bigr)^{-1} 
(1+qW)^{-1}(1+qW + q^2WY + q^3WYZ)
$$
$$
W'= 
W^{-1}(1+qZ)^{-1}(1+qZ + q^2ZX + q^3ZXW)
WY(1+qZ)\Bigl(1+q^{-1}W+q^{-2}YW +q^{-3} ZYW\Bigr)^{-1}
$$
\vskip 3mm

\section{Appendix: the configuration space of $5$ flags in ${\mathbb P}^2$ is of cluster type $E_7$.}

If mutating  the $X$-coordinates  
(as in  (\ref{5.18.04.2as})) 
we get only a finite  collection of different coordinate systems, 
we say that the corresponding cluster ${\cal X}$-space is 
{\it 
of finite type}. Based on  the Fomin-Zelevinsky 
classification theorem \cite{FZII}, 
we showed in \cite{Cluster} that cluster ${\cal X}$-spaces 
of finite type are also parametrised by Dynkin diagrams 
of the Cartan-Killing type.  In this paper, we defined  
 cluster ${\cal X}$-spaces 
under a simplifying assumption that the matrix 
$\varepsilon_{ij}$ is skew-symmetric. In the finite type case 
this boils down to the condition that  the Dynkin diagram is simply-laced. 

Every skew-symmetric matrix $\varepsilon:J\times J \rightarrow \mathbb Z$ with integer entries determines a graph with oriented edges. Namely take $J$ as the set of vertices and connect vertices $i$ and $j$ by $|\varepsilon^{ij}|$ edges oriented towards $j$ if $\varepsilon_{ij}>0$ and towards $i$ otherwise. Conversely any graph with oriented edges and such that edges connecting any two vertices have the same orientation we can associate a skew-symmetric matrix. Graph notation for skew symmetric matrices is sometimes more convenient than the matrix notation. Below we  use this to 
picture a function $\varepsilon_{ij}$ from 
a seed $\mathbf I  = (I, \varepsilon)$ by an oriented graph. 
A mutation is determined by a vertex $k$, pictured by a little 
circle around the vertex. 
The mutated function $\varepsilon'_{ij}$ is given by formula 
(\ref{5.18.04.2}). If $\varepsilon_{kl}=0$, the 
mutations at the vertices $k,l$ commute, so we may picture both of them on the same diagram. 
The mutated graph is shown on the right 
of the initial one. Finally, a seed $\mathbf I  = (I, \varepsilon)$ 
is of finite type 
if mutating the graph of $\varepsilon$ 
we can get a graph isomorphic to the Dynkin diagram of 
type $A_n, D_n, E_6, E_7, E_8$, oriented in a certain way. 

\begin{figure}[ht]
\centerline{\epsfbox{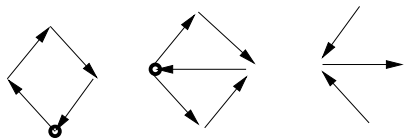}}
\caption{The moduli space of configurations of 
$4$ flags in $PGL_3$ is of cluster 
type $D_4$.}
\label{E7}
\end{figure}


\begin{proposition} \label{finite} a) The moduli space of configurations of $4$ flags in $PGL_3$ 
has a structure of a cluster ${\cal X}$-space of finite type $D_4$.

b). The moduli space of configurations of $5$ flags in $PGL_3$ 
has a structure of a cluster ${\cal X}$-space of finite type $E_7$.
\end{proposition}

\begin{figure}[ht]
\centerline{\epsfbox{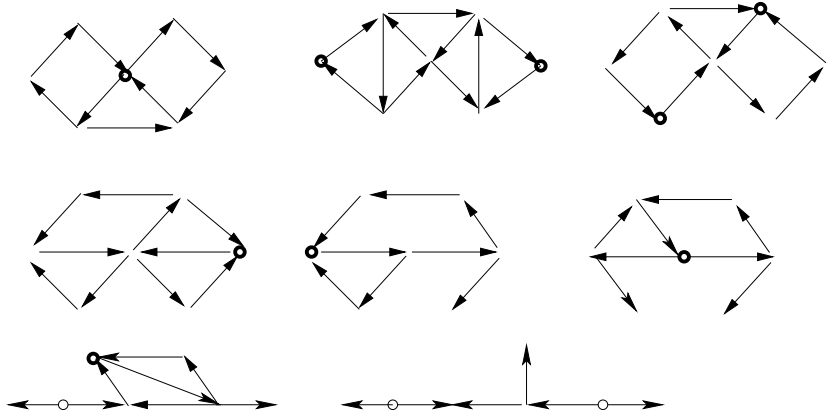}}
\caption{Mutations 
leading to a standard $\varepsilon$-function  of type $E_7$.}
\label{E7}
\end{figure}


{\bf Proof}. We exhibit a sequence 
of mutations of the $\varepsilon$-function  transforming
the $\varepsilon$-function describing the moduli space of configurations 
of $4$ (respectively $5$) flags in $PGL_3$ to a standard one 
of finite type $D_4$ (respectively $E_7$). 
The case $n=4$ is shown on Figure \ref{D4}, and the  $n=5$ case 
on Figure \ref{E7}. 
The proposition is proved.

\vskip 3mm
{\bf Acknowledgments}. This work has been done when 
V.F. visited Brown University. He would like to thank Brown University 
for hospitality and support. He wish to thank grants 
 RFBR 01-01-00549, 0302-17554 and NSh 1999.2003.2. 
A.G. was 
supported by the  NSF grants  DMS-0099390 and DMS-0400449. 
The paper was finished when A.G. enjoyed the hospitality of  
the MPI (Bonn). 
He is grateful to NSF and MPI for the support. 

We are very grateful to the Referee for many useful comments 
 incorporated in the paper.

\end{document}